\documentclass[11pt, letterpaper]{article}
\usepackage[utf8]{inputenc}
\usepackage{microtype}
\usepackage{graphicx}
\usepackage{subfigure}
\usepackage{booktabs}
\usepackage{hyperref}
\usepackage{amsmath}
\usepackage{amssymb}
\usepackage{color}
\usepackage{authblk}
\usepackage[margin=1in]{geometry}
\usepackage{cite}

\usepackage{tikz}
\usetikzlibrary{shapes,arrows}
\usetikzlibrary{positioning}
\tikzstyle{block} = [draw, fill=white, rectangle, 
    minimum height=3em, minimum width=3em]
\tikzstyle{circ} = [draw, fill=white, circle, minimum size=2.5em]
\tikzstyle{input} = [coordinate]
\tikzstyle{output} = [coordinate]
\tikzstyle{pinstyle} = [pin edge={to-,thin,black}]

\title{Connections Between Adaptive Control and Optimization in Machine Learning}
\author[1]{Joseph E. Gaudio}
\author[2]{Travis E. Gibson}
\author[1]{Anuradha M. Annaswamy}
\author[3]{Michael A. Bolender}
\author[4]{Eugene Lavretsky}
\affil[1]{Massachusetts Institute of Technology}
\affil[2]{Brigham and Women’s Hospital and Harvard Medical School}
\affil[3]{Air Force Research Laboratory}
\affil[4]{The Boeing Company}
\date{April 11, 2019}

\begin{document}

\maketitle

\begin{abstract}
This paper demonstrates many immediate connections between adaptive control and optimization methods commonly employed in machine learning. Starting from common output error formulations, similarities in update law modifications are examined. Concepts in stability, performance, and learning, common to both fields are then discussed. Building on the similarities in update laws and common concepts, new intersections and opportunities for improved algorithm analysis are provided. In particular, a specific problem related to higher order learning is solved through insights obtained from these intersections.
\end{abstract}

\section{Introduction}
\label{s:INTRODUCTION}

The fields of adaptive control and machine learning have evolved in parallel over the past few decades, with a significant overlap in goals, problem statements, and tools. Machine learning as a field has focused on computer based systems that improve through experience \cite{Duda_2001,Bishop_2006,Hastie_2009,Efron_2016,Goodfellow-et-al-2016,Jordan_2015}. Often times the process of learning is encapsulated in the form of a parameterized model, whose parameters are learned in order to approximate a function. Optimization methods are commonly employed to reduce the function approximation error using any and all available data. The field of adaptive control, on the other hand, has focused on the process of controlling engineering systems in order to accomplish regulation and tracking of critical variables of interest (e.g. speed in automotive systems, position and force in robotics, Mach number and altitude in aerospace systems, frequency and voltage in power systems) in the presence of uncertainties in the underlying system models, changes in the environment, and unforeseen variations in the overall infrastructure \cite{Narendra_1989,Sastry_1989,Astrom_1995,Ioannou1996,Narendra2005}. The approach used for accomplishing such regulation and tracking in adaptive control is the learning of underlying parameters through an online estimation algorithm. Stability theory is employed for enabling guarantees for the safe evolution of the critical variables, and convergence of the regulation and tracking errors to zero.

Learning parameters of a model in both machine learning and adaptive control occurs through the use of input-output data. In both cases, the main algorithm used for updating the parameters is based on a gradient descent-like algorithm \cite{Narendra2005}. Related tools of analysis, convergence, and robustness in both fields have a tremendous amount of similarity. As the scope of problems in both fields increases, the associated complexity and challenges increase as well. Therefore it is highly attractive to understand these similarities and connections so that the two communities can develop new methods for addressing new challenges. In this paper, we discuss the similarities and connections in detail between the fields of adaptive control and machine learning. Using these connections, we state and provide a solution for a new problem in machine learning using methods developed in adaptive control.

In this paper, the adaptive control perspective will be presented in continuous time with machine learning material presented in discrete time. The paper organization is as follows. We introduce the formulation of output errors commonly employed in adaptive control and machine learning with their associated update laws in Section \ref{s:SETUP}. Numerous connections between the two fields are then made with respect to the underlying parameter update laws in Section \ref{s:CONNECTIONS_Update_Law} and important concepts in Section \ref{s:CONNECTIONS_Concepts}. Examples of intersections between both fields are provided in Section \ref{s:WHY}, with concluding remarks in Section \ref{s:CONCLUSIONS}.

\section{Problem Statements}
\label{s:SETUP}

In this section, we state typical problems that are addressed in the areas of adaptive control and machine learning. In both cases, we illustrate the role of learning, the input-output data used, and the overall problem that is desired to be solved.

\subsection{Adaptive Control}
\label{ss:SETUP_Adaptive_Control}

The main goal in adaptive control is to carry out problems such as estimation or tracking in the presence of parametric uncertainties. The underlying model that relates inputs, outputs, and the unknown parameters is assumed to stem from either the underlying physics or from data-driven approaches. Often these models take the form
\begin{equation}\label{e:y_algebraic}
    y(t)=f_1(\phi(t),\theta^*)
\end{equation}
or
\begin{equation}\label{e:y_dynamical}
    \dot{x}(t)=f_2(x(t),u(t),\theta^*),\quad y(t)=f_3(x(t),u(t),\theta^*)
\end{equation}
where $u\in\mathbb{R}^m$ is an exogenous input, $x\in\mathbb{R}^n$ denotes the state, $y\in\mathbb{R}^p$ corresponds to output measurements, $\phi\in\mathbb{R}^N$ corresponds to measured and computed variables, and $\theta^*\in\mathbb{R}^N$ denotes the uncertain parameter. In an estimation problem, the goal is to estimate the state $x$ in (\ref{e:y_dynamical}) and output $y$ in both (\ref{e:y_algebraic}), (\ref{e:y_dynamical}), alongside the unknown parameter $\theta^*$ simultaneously, using all available variables. In a control problem, the goal is to determine a control input $u$ so that the output $y$ in (\ref{e:y_dynamical}) follows a desired output $\hat{y}$. 

A typical approach taken in order to solve the estimation problem in (\ref{e:y_algebraic}) is to choose an estimator structure of the form
\begin{equation}\label{e:y_hat_algebraic_AC}
    \hat{y}(t)=f(\phi(t),\theta(t)) 
\end{equation}
where $\theta\in\mathbb{R}^N$ denotes the estimate of $\theta^*$ and adjust $\theta$ so that the estimation error $e_y=\hat{y}-y$ is minimized, i.e., choose a function $g_1(e_y,\phi)$ with
\begin{equation}\label{e:g_1}
    \dot{\theta}(t)=g_1(e_y(t),\phi(t))
\end{equation}
so that the estimator has bounded signals, $e_y(t)$ converges to zero and $\theta(t)$ converges to $\theta^*$. Similarly, the control problem consists of constructing an output tracking error $e_y=\hat{y}-y$, where $\hat{y}$ denotes the desired output that $y$ is required to track. The goal is to then choose functions $g_2(e_y,\phi,\theta)$ and $g_3(e_y,\phi,\theta)$ so that the control input $u$ and a control parameter estimate $\theta$ can be chosen as
\begin{align}
    \begin{split}
        u(t)&=g_2(e_y(t),\phi(t),\theta(t))\\
        \dot{\theta}(t)&=g_3(e_y(t),\phi(t),\theta(t))
    \end{split}
\end{align}
leading to closed-loop signals remaining bounded, $e_y(t)$ converging to zero and $\theta(t)$ converging to its true value $\theta^*$. Denote the corresponding parameter errors as $\tilde{\theta}=\theta-\theta^*$.

In order to derive the function $g_1$ for the estimation problem in (\ref{e:y_algebraic}) and the functions $g_2$ and $g_3$ for the control problem in (\ref{e:y_dynamical}) so as to realize the underlying goals, a stability framework together with an error model approach is often employed in adaptive control. The error model approach consists of identifying the basic relationship between the two errors that are commonly present in these adaptive systems, which are the estimation (or tracking) error $e_y$ and the parameter error $\tilde{\theta}$. While the estimation error is measurable and correlated with the parameter error, the parameter error is unknown but adjustable through the parameter estimate. In order to determine the update laws $g_i$, the relationship (error model) that relates these two errors is used as a cue. 

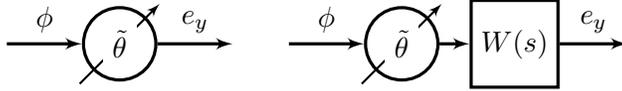
\begin{figure}[t]
    \vspace{0.15cm}
	\centering
    \begin{tikzpicture}[auto, node distance=1.5cm,>=latex']
    
    \node [input, name=input1] {};
    \node [circ, right of=input1, very thick] (circle1) {};
    \node [output, name=output1, right of=circle1] {};
    \draw [->, very thick] (input1) -- node [name=phi1] {$\phi$} (circle1);
    \draw [->, very thick] (circle1) -- node [name=e_y1] {$e_y$} (output1);
    \node [input, name=arrowleft1, below left=0.25cm of circle1] {};
    \node [input, name=arrowright1, above right=0.25cm of circle1] {};
    \draw [->, thick] (arrowleft1) -- (arrowright1);
    \node[fill=white] at (circle1) {$\tilde{\theta}$};
    
    \node [input, name=input2, right of=output1, node distance=0.75cm] {};
    \node [circ, right of=input2, very thick] (circle2) {};
    \node [block, right of=circle2, very thick] (plant) {$W(s)$};
    \node [output, name=output2, right of=plant] {};
    \draw [->, very thick] (input2) -- node [name=phi2] {$\phi$} (circle2);
    \draw [->, very thick] (circle2) -- (plant);
    \draw [->, very thick] (plant) -- node [name=e_y_2] {$e_y$} (output2);
    \node [input, name=arrowleft2, below left=0.25cm of circle2] {};
    \node [input, name=arrowright2, above right=0.25cm of circle2] {};
    \draw [->, thick] (arrowleft2) -- (arrowright2);
    \node[fill=white] at (circle2) {$\tilde{\theta}$};
    
    \end{tikzpicture}
	\caption{Error Models. \textbf{Left:} Regression (\ref{e:error_model_1}), \textbf{Right:} Adaptive Control (\ref{e:error_model_3}).}
	\label{f:Block_Diagram_Error_1_and_3}
\end{figure}

Two types of error models frequently occur in adaptive systems, and are presented below (see Figure \ref{f:Block_Diagram_Error_1_and_3}). The first corresponds to the case when the relation in (\ref{e:y_algebraic}) is linear, and the underlying error model is simply of the form (cf. \cite{Narendra2005})
\begin{equation}\label{e:error_model_1}
    e_y(t)=\tilde{\theta}^T(t)\phi(t)
\end{equation}
and as a result, the function $g_1$ in (\ref{e:g_1}) can be determined simply using the gradient rule that minimizes $\lVert e_y\rVert^2$. The second is of the form (cf. \cite{Narendra2005})
\begin{equation}\label{e:error_model_3}
    e_y(t)=W(s)[\tilde{\theta}^T(t)\phi(t)]
\end{equation}
where $W(s)[\zeta]$ denotes a dynamic operator operating on $\zeta(t)$. It has been shown in the adaptive control literature \cite{Narendra_1989,Sastry_1989,Astrom_1995,Ioannou1996,Narendra2005} that for specific classes of dynamic operators $W(s)$, a stable, gradient-like rule can be determined for adjusting $\tilde{\theta}$. Most of these results apply uniformly to the case when $u$ and $y$ are scalars or vectors, with the latter introducing additional technicalities. In this paper we consider the case where inputs and outputs are scalars for notational simplicity, and to focus on the core of the learning problem with multi-dimensional regressors $\phi$ and parameter estimates $\theta$. Often the unknown parameter $\theta^*$ is assumed to reside in a compact convex set, which we will denote as $\Theta$.

\subsection{Machine Learning}
\label{ss:SETUP_Machine_Learning}

Machine learning is a broad field encompassing a wide variety of learning techniques and problems such as classification and regression. A large portion of machine learning considers supervised learning problems, where regressors $\phi$ and outputs $y$ are related to one another in an unknown algebraic manner \cite{Duda_2001,Bishop_2006,Hastie_2009,Efron_2016,Goodfellow-et-al-2016,Jordan_2015}. A typical approach taken in order to perform classification or regression is to choose an output estimator $\hat{y}_k$ parameterized with adjustable weights $\theta_k$ as
\begin{equation}\label{e:y_hat_ML}
    \hat{y}_k=f(\phi_k,\theta_k).
\end{equation}
A common form of the estimator as in (\ref{e:y_hat_ML}) is that of neural networks, where the parameters $\theta_k$ represent the adjustable weights in the network \cite{Duda_2001,Bishop_2006,Hastie_2009,Efron_2016,Goodfellow-et-al-2016}.

Similar to adaptive control, $\theta_k$ is often adjusted using the output error $e_{y,k}=\hat{y}_k-y_k$. A loss function $L:\Theta\rightarrow\mathbb{R}$ of $e_{y,k}$ is minimized through the adjustable weights. An example loss function for regression is $\ell_p$ loss (with $p\in\mathbb{N}$, $p>0$ and even) $L(\theta_k)=(1/p)\lVert e_{y,k}\rVert_p^p$. For binary classification ($y_k\in\{-1,1\}$) common loss functions include hinge loss $L(\theta_k)=\max(0,1-y_k\hat{y}_k)$, and logistic loss $L(\theta_k)=\ln(1+\text{exp}(-y_k\hat{y}_k))$. Additionally, as in empirical risk minimization (ERM) \cite{Vapnik_1992}, the total loss function considered for the purpose of a parameter update may be an average of loss functions over $m$ samples as: $(1/m)\sum_{i=1}^mL_i(\theta_k)$. The above descriptions make it clear that the structure of the estimation problem in both adaptive control and machine learning are strikingly similar. In the next section, we examine the nature of the adjustment of $\theta_k$.

\subsection{Common Update Laws}
\label{ss:SETUP_Update_Laws}

As previously stated, the goal in adaptive control is to design a rule to adjust $\theta$ in an online continuous manner using knowledge of $\phi$ and $e_y$ such that $e_y$ tends toward zero. Given that the output errors may be corrupted by noise, an iterative, gradient-like update is usually employed. To do so for the algebraic error model (\ref{e:error_model_1}), consider the squared loss cost function: $L(\theta(t))=(1/2)e_y^2(t)$. The gradient of this function with respect to the parameters can be expressed as: $\nabla_{\theta}L(\theta(t))=\phi(t)e_y(t)$. The standard gradient flow update law \cite{Narendra_1989} may be expressed as follows with user-designed gain parameter $\gamma>0$ as
\begin{equation}\label{e:update_GF}
\dot{\theta}(t)=-\gamma\nabla_{\theta}L(\theta(t))=-\gamma \phi(t)e_y(t).
\end{equation}

For dynamical error models such as (\ref{e:error_model_3}), a stability approach rather than a gradient based one is taken using Lyapunov methods, which leads to an adaptive law identical to (\ref{e:update_GF}) for a class of dynamic systems $W(s)$ that are strictly positive real \cite{Narendra_1989,Parks_1966}.

The common update law for supervised machine learning problems, gradient descent\footnote{While this is not true of all machine learning as the field is broad, (for example Bayesian methods often use sampling based techniques such as Markov Chain Monte Carlo), even in the world of probabilistic inference, gradient based methods can also be used, cf. variational inference \cite{Blei_2017}.}, is akin to the time varying regression law (\ref{e:update_GF}) in discrete time, and of the form
\begin{equation}\label{e:gradient_descent}
    \theta_{k+1}=\theta_k-\gamma_k\nabla_{\theta} L(\theta_k)
\end{equation}
where the ``stepsize" $\gamma_k$ is usually chosen as a decreasing function of time \cite{Hazan_2007,Hazan_2008,Hazan_2016,Bubeck_2015,Zinkevich_2003}, a standard feature of stochastic gradient algorithms.

\section{Connections: Update Law}
\label{s:CONNECTIONS_Update_Law}

This section details a variety of connections between adaptive control and the optimization methods commonly used in machine learning as viewed from the perspective of their common update laws (\ref{e:update_GF}), (\ref{e:gradient_descent}).

\subsection{$\sigma$-Modification, $e$-Modification, and Regularization}

While the update laws in (\ref{e:update_GF}) and (\ref{e:gradient_descent}) are designed primarily to reduce the output error $e_y$, there are several secondary reasons to modify these update laws from robustness considerations due to perturbations stemming from disturbances, noise, and other unmodeled causes. We outline these updates in in this section.

\subsubsection{Adaptive Control}

Historically the adaptive update law in (\ref{e:update_GF}) has been modified to ensure robustness in the presence of bounded disturbances as
\begin{equation}\label{e:sigma_e_mod}
\dot{\theta}(t)=-\gamma\left[\nabla_{\theta}L(\theta(t))+\sigma \mathcal{G}(\theta(t),e_y(t))\right]
\end{equation}
where $\sigma>0$ is a tuneable parameter that scales the extra term $\mathcal{G}$. Common choices for $\mathcal{G}$ include the $\sigma$-modification $\mathcal{G}=\theta$ \cite{Ioannou_1984}, and the $e$-modification $\mathcal{G}=\lVert e_y\rVert\theta$ \cite{Narendra_1987a}.

\subsubsection{Machine Learning}

Regularization is often included in a machine learning optimization problem in order to help cope with overfitting by including constraints on the parameter, thus resulting in an augmented loss function \cite{Duda_2001,Bishop_2006,Efron_2016,Goodfellow-et-al-2016,Hazan_2008,Hastie_2009,Bubeck_2015,Hazan_2016}: $\bar{L}(\theta)=L(\theta)+\sigma\mathcal{R}(\theta)$ where $\sigma>0$ is a tunable parameter, often referred to as a Lagrange multiplier. The gradient descent update (\ref{e:gradient_descent}) for this augmented loss function is often referred to as the ``regularized follow the leader" algorithm in online learning \cite{Hazan_2016} and may be expressed as
\begin{equation}
    \theta_{k+1}=\theta_k-\gamma_k\left[\nabla_{\theta} L(\theta_k)+\sigma\nabla_{\theta}\mathcal{R}(\theta_k)\right].
\end{equation}
The common choice of $\ell_p$ regularization in machine learning of $\mathcal{R}=(1/p)\lVert\theta\rVert_p^p$ with $p=2$, (as in ridge regression), coincides with the $\sigma$-modification \cite{Ioannou_1984}, as then $\nabla_{\theta}\mathcal{R}=\mathcal{G}$. Given that the dimension of the parameter vector may be large, a sparse representation is often obtained with $\ell_1$ regularization (as in lasso), with $\mathcal{R}=\lVert\theta\rVert_1$ \cite{Bishop_2006,Hastie_2009,Efron_2016,Goodfellow-et-al-2016}.

\subsection{Deadzone Modification and Early Stopping}

This subsection details common modifications of the adaptive law adopted to cease updating the parameter estimate after sufficient tuning.

\subsubsection{Adaptive Control}

Another method employed to increase robustness in the presence of bounded disturbances is to employ a ``dead zone" \cite{Peterson_1982}, for the update in (\ref{e:update_GF}) as
\begin{equation}
    \dot{\theta}(t)=
    \left\{\begin{array}{ll}
        \hspace{-.15cm}-\gamma\nabla_{\theta}L(\theta(t)), & \mathcal{D}(e_y)>d_0+\epsilon\\
        \hspace{-.15cm}0, & \mathcal{D}(e_y)\leq d_0+\epsilon
        \end{array}\right.
\end{equation}
where $d_0>0$ is the dead zone width that may correspond to an upper bound on the disturbance, and $\epsilon>0$ is a small constant. The function $\mathcal{D}$ is a non-negative metric on the output error to stop adaptation in desired regions of the output space. A common choice is $\mathcal{D}=\lVert e_y\rVert$ such that adaptation stops after a small output error is achieved above a noise level with upper bound $d_0$.

\subsubsection{Machine Learning}

The training processes is often stopped in machine learning applications as a method to deal with overfitting \cite{Bishop_2006,Hastie_2009,Efron_2016,Goodfellow-et-al-2016,Prechelt_1998}. This may be done by using multiple data sets and stopping the parameter update process (\ref{e:gradient_descent}) when the loss computed for a validation data set starts to increase \cite{Prechelt_1998}. Early stopping is often seen to be needed for training neural networks due to their large number of parameters \cite{Bishop_2006,Hastie_2009,Efron_2016,Goodfellow-et-al-2016} and can act as regularization \cite{SJ_BERG_1995}.

\subsection{Projection}

It is often desirable to define a compact region a priori for the parameters $\theta$, such that during the learning process the parameters are not allowed to leave that region. In physical systems there are natural constraints which may aid in the design of that region, and for non physical systems, the constraints are often engineered by the algorithm designer.

\subsubsection{Adaptive Control}
\label{sss:Projection_Adaptive_Control}

A continuous projection algorithm is commonly employed to provide for robustness of the adaptive update law in the presence of unmodeled dynamics \cite{Kreisselmeier_1982,Hussain_2017,Gibson_2012}. One such implementation is
\begin{equation}\label{e:Projection}
    \text{Proj}(\theta_i,\zeta_i)=
    \left\{\begin{array}{ll}
        \hspace{-.15cm}\frac{\theta_{i,\max}^2-\theta_i^2}{\theta_{i,\max}^2-\theta_{i,\max}'^2}\zeta_i, & \hspace{-0.11cm}\theta_i\in\Omega_i\wedge\theta_i\zeta_i>0\\
        \hspace{-.15cm}\zeta_i, & \hspace{-0.11cm}\text{otherwise}
        \end{array}\right.
\end{equation}
where $\Omega$, $\theta_{i,\max}$, $\theta_{i,\max}'$ define a user-specified boundary layer region inside of $\Theta$ (see \cite{Hussain_2017}). The update law in (\ref{e:update_GF}) may then be modified as
\begin{equation}
    \dot{\theta}(t)=-\gamma\text{Proj}\left[\theta(t),\nabla_{\theta}L(\theta(t))\right].
\end{equation}

\subsubsection{Machine Learning}

Projected gradient descent methods have a long history in optimization. The following projection operation finds the point in a convex set which is closest to a specified point, and may be defined as
\begin{equation}
    \Pi_{\Theta}(\bar{\theta})\triangleq\underset{\theta\in\Theta}{\arg\min}\lVert\theta-\bar{\theta}\rVert
\end{equation}
which may be employed in the update sequence \cite{Hazan_2007,Hazan_2008,Hazan_2016,Bubeck_2015,Zinkevich_2003}
\begin{equation}\label{e:projected_GD}
    \bar{\theta}_{k+1}=\theta_k-\gamma_k\nabla_{\theta}L(\theta_k),\qquad\theta_{k+1}=\Pi_{\Theta}(\bar{\theta}_{k+1}).
\end{equation}

\subsection{Adaptive Gains and Stepsizes}

\subsubsection{Adaptive Control}

The following parameter update law for the algebraic error model\footnote{This update law has not been proven stable for the error model in (\ref{e:error_model_3}).} (\ref{e:error_model_1}) is one example which alters the gain of the standard update law (\ref{e:update_GF}) as a function of the time varying regressors $\phi$ \cite{Narendra_1989,Ioannou1996}:
\begin{align}\label{e:time_varying_gain}
    \begin{split}
        \dot{\theta}(t)&=-\gamma\Gamma(t)\nabla_{\theta}L(\theta(t))\\
        \dot{\Gamma}(t)&=
        \left\{\begin{array}{ll}
        \hspace{-0.195cm}\Upsilon \Gamma(t) -\frac{\Gamma(t)\phi(t)\phi^T(t)\Gamma(t)}{\mathcal{N}(t)}, & \hspace{-0.11cm}\lVert\Gamma(t)\rVert\leq\Gamma_{\max}\\
        \hspace{-0.16cm}0, & \hspace{-0.11cm}\text{otherwise}
        \end{array}\right.
    \end{split}
\end{align}
where $\Upsilon\geq0$ is a \emph{forgetting factor} and $\mathcal{N}(t)$ is a \emph{normalizing signal}, with common choice $\mathcal{N}(t)=(1+\mu\phi^T(t)\phi(t))$ for $\mu>0$ chosen appropriately (see for example \cite{Ioannou1996} for a discussion of the choice of parameters). It can be seen that the update for $\Gamma$ may be integrated and used in the update for $\theta$ to result in a gain adaptive to the regressor $\phi$.

\subsubsection{Machine Learning}

Adaptive step size methods \cite{Duchi_2011,Zeiler_2012,Kingma_2017,Reddi_2018} have seen widespread use in machine learning problems due to their ability to handle sparse and small gradients by adjusting the step size as a function of features as they are processed online. Define the following: $g_k=\nabla_{\theta} L(\theta_k)$, $m_k=\mathcal{F}_{1,k}(g_1,\ldots,g_k)$, $V_k=\mathcal{F}_{2,k}(g_1,\ldots,g_k)$ for user defined averaging functions $\mathcal{F}_{1,k}$, $\mathcal{F}_{2,k}$. A common update law for adaptive step size methods \cite{Reddi_2018} can then be seen to be similar to (\ref{e:projected_GD}) as
\begin{equation}\label{e:adaptive_stepsize_generic}
    \bar{\theta}_{k+1}=\theta_k-\gamma_km_k/V_k^{1/2},\qquad\theta_{k+1}=\Pi_{\Theta}(\bar{\theta}_{k+1})
\end{equation}
where the following parameterizations are common \cite{Reddi_2018}: (i) projected gradient descent\footnote{$\mathcal{F}_{1,k}=g_k$, $\mathcal{F}_{2,k}=I$.} (\ref{e:projected_GD}), (ii) \textsc{AdaGrad}\footnote{$\mathcal{F}_{1,k}=g_k$, $\mathcal{F}_{2,k}=\epsilon I+\text{diag}(\sum_{i=1}^kg_i^2)$, where $g_i^2=g_i\odot g_i$.} \cite{Duchi_2011}, and (iii) \textsc{Adam}\footnote{$\mathcal{F}_{1,k}=(1-\beta_1)\sum_{i=1}^k\beta_1^{k-i}g_i$, $\mathcal{F}_{2,k}=(1-\beta_2)\text{diag}(\sum_{i=1}^k\beta_2^{k-i}g_i^2)$.} \cite{Kingma_2017}. It can be noted that the normalization in these update laws is a function of the gradient, which can be compared to the normalization by the regressor in (\ref{e:time_varying_gain}).

\section{Connections: Tools and Concepts}
\label{s:CONNECTIONS_Concepts}

This section details concepts and tools common to both machine learning and adaptive control.

\subsection{Lyapunov Functions and Regret}

Stability and convergence tools in adaptive control and online machine learning are analyzed in this section.

\subsubsection{Adaptive Control}

Suppose we consider the error model in (\ref{e:error_model_3}) where $W(s)=c(sI-A)^{-1}b$, and a corresponding state space representation of the form \cite{Narendra_1989}
\begin{align}
\begin{split}
\label{e:error_model_3_general}
\dot{e}(t)&=Ae(t)+b\tilde{\theta}^T(t)\hat{\phi}(t)+b\theta^{*T}\tilde{\phi}(t)\\
e_y(t)&=ce(t).
\end{split}
\end{align}
The term $\tilde{\phi}$ is due to exponentially decaying terms in the regressor $\phi$. That is, $\tilde{\phi}=\hat{\phi}-\phi$ and $\dot{\tilde{\phi}}=\Lambda\tilde{\phi}$ for a Hurwitz matrix $\Lambda\in\mathbb{R}^{N\times N}$.\footnote{This formulation is common in the design of non-minimal adaptive observers \cite{Narendra_1989}. It can be noted that $\hat{\phi}\rightarrow\phi$ as $t\rightarrow\infty$ as $\Lambda$ is Hurwitz. Also for $\hat{\phi}=\phi$, (\ref{e:error_model_3_general}) is the same as (\ref{e:error_model_3}). A Hurwitz matrix $\Lambda$ implies the existence of a positive definite matrix $\bar{P}=\bar{P}^T\in\mathbb{R}^{N\times N}$ and $0<\bar{Q}=\bar{Q}^T\in\mathbb{R}^{N\times N}$ such that: $\Lambda^T\bar{P}+\bar{P}\Lambda=-\bar{Q}$.} Stability is often proven in adaptive control by the use of a Lyapunov function $V$, such as
\begin{equation}
    V=\gamma^{-1}\tilde{\theta}^T\tilde{\theta}+e^TPe+\alpha\tilde{\phi}^T\bar{P}\tilde{\phi}.
\end{equation}
It should be noted that the last two terms in $V$ are not needed for the algebraic error model in (\ref{e:error_model_1}). The time derivative of the Lyapunov function may then be stated using the update law in (\ref{e:update_GF}) and the KYP lemma \cite{Narendra_1989} as
\begin{equation}\label{e:v_dot}
    \dot{V}=-e^TQe-\alpha\tilde{\phi}^T\bar{Q}\tilde{\phi}+2e^TPb\theta^{*T}\tilde{\phi}
\end{equation}
where $\dot{V}\leq0$ for $\alpha>(4\lVert Pb\rVert^2\lVert\theta^*\rVert^2/(\min eig(Q)\cdot\min eig(\bar{Q}))$. It can be shown \cite{Narendra_1989} that $\delta(t)\triangleq 2e^TPb\theta^{*T}\tilde{\phi}$ is an exponentially decaying signal with $\tilde{\phi},e\in\mathcal{L}_2\cap\mathcal{L}_{\infty}$. By integrating $\dot{V}$ from $t_0$ to $T$, we obtain
\begin{equation}\label{e:regret_AC}
    \int_{t_0}^{T}e^TQedt-\int_{t_0}^{T}\delta(t)dt\leq-\int_{t_0}^{T}\dot{V}dt=V(t_0)-V(T).
\end{equation}
Given that $\dot{V}\leq0$, $V(t_0)-V(T)\leq V(t_0)<\infty$.

\subsubsection{Machine Learning}\label{sss:Regret}

In online learning, efficiency of an algorithm is often analyzed using the notion of ``regret" as
\begin{equation}\label{e:regret}
    \text{regret}_T=\sum_{k=1}^T\mathcal{C}_k(\theta_k)-\min_{\theta\in\Theta}\sum_{k=1}^T\mathcal{C}_k(\theta)
\end{equation}
where regret can be seen to correspond to the sum of the time varying convex costs $\mathcal{C}_k$ associated with the choice of the time varying parameter estimate $\theta_k$, minus the cost associated with the best static parameter estimate choice, over a time horizon of $T$ steps \cite{Hazan_2007,Hazan_2008,Hazan_2016,Zinkevich_2003}. Suppose we consider a quadratic cost $\mathcal{C}_k=e_k^TQe_k$, $Q=Q^T>0$. A continuous time limit of (\ref{e:regret}) leads to an integral as
\begin{equation}\label{e:regret_continuous}
    \text{continuous regret}_T=\int_{t_0}^T e^TQedt-\int_{t_0}^T\bar{\delta}(t)dt
\end{equation}
where $\bar{\delta}(t)$ is an exponentially decaying signal which is due to nonzero initial conditions in (\ref{e:error_model_3}) or similarly in (\ref{e:error_model_3_general}).\footnote{This may be seen by setting $\theta(t)\equiv\theta^*$ in (\ref{e:error_model_3}) or (\ref{e:error_model_3_general}), thus resulting in an exponentially decreasing $e^TQe$. Note that this exponentially decaying term is absent in the time varying regression case (\ref{e:error_model_1}).} A strong similarity can thus be seen between (\ref{e:regret_AC}) and (\ref{e:regret_continuous}).

It is desired to have regret grow sub-linearly with time, such that average regret, $(1/T)\text{regret}_T$, goes to zero in the limit $T\rightarrow\infty$, to provide for an efficient algorithm \cite{Hazan_2016}. Average regret can be connected to convergence in the case of a constant cost and by applying Jensen's inequality as \cite{Hazan_2016}
\begin{equation}\label{e:cost_convergence}
    \mathcal{C}(\bar{\theta}_T)-\mathcal{C}(\theta^*)\leq\frac{1}{T}\sum_{k=1}^T\left[\mathcal{C}(\theta_k)-\mathcal{C}(\theta^*)\right]=\frac{\text{regret}_T}{T}
\end{equation}
where $\bar{\theta}_T=(1/T)\sum_{k=1}^T\theta_k$ is the average parameter estimate. Here sub-linear regret helps show convergence of the costs in (\ref{e:cost_convergence}). For adaptive control, convergence of state/output errors is shown from a similar integral which is akin to \emph{constant} regret upper bounded by $V(t_0)$ in (\ref{e:regret_AC}).

\subsection{Unmodeled Dynamics and Generalization}

This section discusses robustness to unforeseen perturbations such as unmodeled dynamics and unseen data.

\subsubsection{Adaptive Control}

Models used to design adaptive controllers, including the examples of (\ref{e:error_model_1}) and (\ref{e:error_model_3}), are linearized approximations with a certain amount of modeling errors. As such, they may only hold about an operating point and need to contend with unmodeled dynamics. This implies that any stabilizing controllers must be designed to not only adapt to parametric uncertainties, but also be robust to unmodeled dynamics. In addition, constraints on the state and input may also be present in adaptive control problems \cite{Karason_1994,Annaswamy_1995}. Analysis becomes more complicated when considering such unmodeled dynamics and constraints, resulting in non-global guarantees. Many of the update law modification in adaptive control from Section \ref{s:CONNECTIONS_Update_Law} were initially derived to ensure robustness in such cases.

\subsubsection{Machine Learning}

This same notion of robustness to modeling errors exists in machine learning in which an estimator $\hat{y}$ is constructed from a finite training data set, often with a finite number of tuneable parameters. It is then desired that this estimator produces a low prediction error based on a test data set consisting of not just seen data, but unseen data as well. Generalization in machine learning thus refers to the concept of a designed estimator having low loss when applied to new problems. In particular it can be seen that in specific cases, generalization pertains to stability, where algorithms that are stable and train in a small amount of time result in a small generalization error \cite{Bousquet_2002,Hardt_2016}.

\subsection{Persistence of Excitation and Stochastic Perturbations}

This section discusses conditions under which parameter estimates can be guaranteed to converge to their true values.

\subsubsection{Adaptive Control}\label{sss:PE_AC}

Persistence of excitation (PE) of the system regressor in adaptive control is a condition that has been shown to be necessary and sufficient for parameter convergence \cite{Jenkins_2018}. It can be shown that if the regressor $\phi$ is persistently exciting, then the algebraic error model (\ref{e:error_model_1}) parameter estimation error $\tilde{\theta}(t)$ converges to zero uniformly in time \cite{Narendra_1989}. Similar conditions can be imposed for the dynamical error model (\ref{e:error_model_3}) and update law (\ref{e:update_GF}) \cite{Narendra_1989}. The PE condition essentially corresponds to certain spectral conditions being satisfied by the regressor \cite{Boyd_1986}.\footnote{In particular, \cite{Boyd_1983} established a condition on spectral lines of signals.} Parameter convergence can also occur through the use of ``the hybrid algorithm", ``the integral algorithm", ``the algorithm with time-varying adaptive gains", and ``the algorithm using multiple models" as is discussed in \cite{Narendra_1989}. A detailed exposition of system identification and parameter convergence in both deterministic and stochastic cases can be found in \cite{Goodwin_1984,Anderson_1982,Narendra_1986,Narendra_1987,Ljung_1987}. Another way to think of the PE condition is that it leads to a perfect test error, since it provides for convergence of the parameter error to zero, and therefore zero output/state error once transients decay to zero.

\subsubsection{Machine Learning}

Many machine learning problems consider the case when stochastic perturbations are present. In this context, significant improvements may be possible by leveraging well known concepts in system identification \cite{Ljung_1987}. For example \cite{Dean_2018} purposely includes a Gaussian random input into a dynamical system in order to provide for PE by construction. Such stochastic perturbations can guarantee a PE condition only in the limit, when infinite samples can be obtained. In order to address the realistic case of finite samples, approaches in machine learning algorithms for system identification and control have attempted to obtain performance bounds with probability $1-p_f$ for $p_f\in(0,1)$.\footnote{The performance bound usually scales inversely with $p_f$ as well.} The probability of failure given by the choice of $p_f$ allows for error due to the presence of finite samples.

\subsection{Tracking vs Exploration}

The concept of exploration can be viewed as the opposite of tracking, with the former often employed in machine learning while the latter is one of the main control goals.

\subsubsection{Adaptive Control}

The goal of adaptive control is to adjust the parameter $\theta$ in such a way to minimize the output error $e_y$ in (\ref{e:error_model_1}) and (\ref{e:error_model_3}). It can be seen from the error models in (\ref{e:error_model_1}) and (\ref{e:error_model_3}) with the update in (\ref{e:update_GF}), that as the output error $e_y$ goes to zero, learning becomes less and less, and that it is possible for a large parameter error to remain even with zero output (or tracking) error. That is, in many adaptive control applications, stability and tracking are successfully accomplished even without parameter convergence.

\subsubsection{Machine Learning}

In many machine learning methods, including reinforcement learning, there exist explicit modifications to update laws to promote exploration of the parameter space. These modifications include restarting trajectories with random initial conditions, adding random perturbations to algorithms, and driving the system towards a non-zero error regions \cite{Dean_2018,Sutton_1992,Sutton_2018}. This preference of exploration and learning over stability is motivated by the desire to find optimal parameters of a system. Stability is not always crucial as models are often trained with offline data on a computer, allowing for many iterations without the financial cost of failure present in physical systems (i.e. a nonzero probability of failure $p_f$ is acceptable).

\subsection{Convergence Guarantees}

Notions of convergence guarantees are of importance in both fields, and are discussed here.

\subsubsection{Adaptive Control}

Adaptive control problems are often parameterized in a specific way such that $e_y$ goes to zero asymptotically as in (\ref{e:error_model_1}) and (\ref{e:error_model_3}). Parameter convergence is shown to occur in these cases with a persistence of excitation condition (see Section \ref{sss:PE_AC}). The specific parameterizations in the output space ensure that a global minimum of $e_y=0$ exists and is unique. In the absence of PE, standard adaptive control algorithms converge to one of the many local minima in the parameter space (i.e. $\dot{\tilde{\theta}}\rightarrow0$ but $\tilde{\theta}\neq0$) \cite{Narendra_1989}.

\subsubsection{Machine Learning}

Machine learning has rapidly grown in recent years, as demonstrated by highly popular and well attended conferences such as ICML and NeurIPS, where rigorous proofs of stability are not always the main focus, instead focusing on empirical performance on large scale problems. A notable exception is a body of work that is emerging which consists of optimization-centric problem formulations, and the examination of the loss landscape, where recent results have shown that in certain classes of problems, local minimums are nearly equivalent to global minimums in terms of performance on test data \cite{Choromanska_2015,Ge_2016,Lee_2016}.

\subsection{Neural Networks}

This section discusses neural networks, a topic common to both fields.

\subsubsection{Adaptive Control}

Gradient based methods to solve for estimates of unknown parameters via back propagation, in what would develop into the foundations of neural networks have been used for decades in control, with early examples consisting of finding optimal trajectories \cite{Pontryagin_1961} in flight control \cite{Kelley_1960}, and resource allocation problems \cite{Bryson_1961} (see \cite{Dreyfus_1990} for a brief history). Since then, the use of neural networks in control systems has expanded to include stabilizing nonlinear dynamical systems \cite{Miller_1995}. Design and analysis of stable controllers based on neural networks was taken up by the adaptive control community due to the the similarities of gradient-like update laws used in neural networks and adaptive control. The adaptive control community developed a well established literature for the use of neural networks in nonlinear dynamical systems in the 1990s \cite{Miller_1995,Narendra_1990,Narendra_1991,Yu_1996,Yu_1998}.

\subsubsection{Machine Learning}

The use of neural networks in the machine learning community greatly expanded as of recent due to the increase in computing power available and an increase in applications \cite{Goodfellow-et-al-2016,Krizhevsky_2012,Sutskever_2013}. Recurrent neural networks \cite{Hopfield_1982,Hinton_1983,Hochreiter_1997}, while often similar in structure to nonlinear dynamical systems, have historically been trained in a manner similar to feed-forward neural networks \cite{Rumelhart_1986} using back propagation through time \cite{Werbos_1990}.\footnote{Hebbian learning \cite{Hebb_1949} based approaches have also been considered.} While a theoretical understanding of why deep neural networks work as well as they do for given problems has been lacking, the machine learning community has worked to rigorously analyze sub-classes of deep neural network architectures such as deep linear networks \cite{Arora_2018,Arora_2019}. The update laws employed in training deep neural networks often include selections of modifications of the update laws as discussed in Section \ref{s:CONNECTIONS_Update_Law}. For an overview of the history of neural networks see \cite{Schmidhuber_2015}.

\subsection{Other Parameterization Schemes}

In addition to neural networks as discussed in the previous section, other parameterizations are often considered in adaptive control and machine learning.

\subsubsection{Adaptive Control}\label{sss:Parameterization_AC}

Adaptive control schemes often consider the case where an unknown parameter occurs linearly with respect to a regressor vector $\phi$ and may be related to an output error $e_y$ algebraically (\ref{e:error_model_1}) or dynamically (\ref{e:error_model_3}). Often times the vector $\phi$ is a nonlinear function of the state of the system or reference model in order to approximate a more general nonlinear function $D$ as: $D(x)=\theta^{*T}\phi(x)$ \cite{Sanner_1992}. Common parameterizations for unknown nonlinearities include Gaussian radial basis functions \cite{Sanner_1992}. Another class of parameterizations consist of nonlinearly parameterized uncertainty $D(\theta^*,\phi)$ in dynamical systems, for which there exists stabilizing adaptive control methods \cite{Loh_1999,Cao_2003}.

\subsubsection{Machine Learning}

Parametric methods are common in machine learning as well, and are useful in many regression and classification based tasks \cite{Duda_2001,Bishop_2006,Hastie_2009,Efron_2016,Goodfellow-et-al-2016}. However, Bayesian based approaches are also widespread in areas such as topic models \cite{Blei_2003}, clustering \cite{Teh_2005} and graphical models \cite{Wainwright_2007}. Additionally, new results in high dimensional statistics are increasingly being considered in which the model may be of higher dimension than the sample size \cite{Wainwright_2019}.

\section{Advantageous Combinations of Machine Learning and Adaptive Control Tools}
\label{s:WHY}

Given the enormous number of similarities in problem statements, tools, concepts, and algorithms, it is natural to examine what the benefits are that accrue by combining insights obtained in these two different communities. Two examples of such an exercise is delineated in this section.

\subsection{Higher Order Learning}\label{ss:HOL}

Many of the update laws addressed thus far were first-order in nature, and made use of gradient-like quantities for learning. A question of increasing interest in the ML community is when accelerated learning can occur for higher-order learning methods. Higher order learning methods are commonly used in machine learning practice \cite{Beck_2009,Krizhevsky_2012,Sutskever_2013} as they can provide for a guaranteed bound on a faster rate of convergence. In particular, Nesterov's accelerated method \cite{Nesterov_1983} was able to certify a convergence rate of $O(1/k^2)$ as compared to the standard gradient descent (\ref{e:gradient_descent}) rate of $O(1/k)$ for a class of convex functions. A parameterization of Nesterov's accelerated method may be stated as
\begin{align}
    \label{e:Nesterov}
    \begin{split}
    \theta_{k+1}&=\vartheta_k-\gamma\nabla_{\theta}L(\vartheta_k)\\
    \vartheta_k&=\theta_k+\beta(\theta_k-\theta_{k-1})
    \end{split}
\end{align}
where $\beta>0$ is a design parameter that weighs the effect of past parameters. Continuous time problem formulations have been explored in \cite{Su_2016,Wibisono_2016}, with rate-matching discretizations established in \cite{Wilson_2016,Wilson_2018,Betancourt_2018}. Many of these methods however become inadequate for time varying inputs.

Adaptive update laws which include additional levels of integration appeared in the ``higher order tuners" in \cite{Morse_1992,Evesque_2003}, and take the form
\begin{align}
    \label{e:Accelerated_GF}
    \begin{split}
    \dot{\vartheta}(t)&=-\gamma\nabla_{\theta}L(\theta(t))\\
    \dot{\theta}(t)&=-\beta(\theta(t)-\vartheta(t))\mathcal{N}(t)
    \end{split}
\end{align}
where $\mathcal{N}(t)\triangleq(1+\mu\phi^T(t)\phi(t))$ for a $\mu>0$. This update law can be seen to be the standard first order update (\ref{e:update_GF}) passed through a time varying filter normalized by the regressor. It was shown in \cite{Gaudio_2019} that (\ref{e:Accelerated_GF}) can provide for rates comparable to accelerated methods in machine learning for static features \cite{Wibisono_2016}. In addition, in contrast to (\ref{e:Nesterov}), the update law in (\ref{e:Accelerated_GF}) can be shown to be stable in the presence of time varying regressors as in (\ref{e:error_model_1}) and as well as in adaptive control applications with error model as in (\ref{e:error_model_3}) \cite{Gaudio_2019}. This extension of accelerated methods in machine learning to include time varying and dynamic error models was only possible by leveraging techniques from adaptive control \cite{Gaudio_2019}.

\subsection{Improved Algorithm Performance Bounds}

Regret analysis common in online machine learning (see Section \ref{sss:Regret}) can result in overly conservative bounds for the performance of an algorithm. In particular, in online projected gradient descent (\ref{e:projected_GD}) for regression (\ref{e:error_model_1}) with squared output error cost $\mathcal{C}=(1/2)e_y^2$, regret analysis guarantees $\text{regret}_T=O(\sqrt{T})$ (cf. \cite{Hazan_2016}). For the same regret cost function, one can guarantee $\text{regret}_T=O(1)$ (constant regret)\footnote{For regression as in (\ref{e:error_model_1}), regret contains a sum of non-negative costs and is therefore a non-decreasing function of the time horizon $T$. Thus $O(1)$ regret is the best achievable regret.}, using adaptive control methods.

\section{Conclusions}
\label{s:CONCLUSIONS}

This paper explored many immediate connections between adaptive control and machine learning, both through common update laws as well as common concepts. Adaptive control as a field has focused on mathematical rigor and guaranteed convergence. The rapid advances in machine learning on the other hand have brought about a plethora of new techniques and problems for learning. This paper was written to elucidate the numerous common connections between both fields such that results from both may be leveraged together to solve new problems.

\section{Acknowledgements}

The authors acknowledge Dr. Michael I. Jordan for useful discussions. This work was supported by the Air Force Research Laboratory, Collaborative Research and Development for Innovative Aerospace Leadership (CRDInAL), Thrust 3 - Control Automation and Mechanization grant FA 8650-16-C-2642 and the Boeing Strategic University Initiative.

\bibliography{References}
\bibliographystyle{ieeetr}

\end{document}